\documentclass[11pt,fleqn]{article}
\newcommand{\ol}{\setlength{\itemsep}{0pt.}\begin{enumerate}}
\newcommand{\eol}{\end{enumerate}\setlength{\itemsep}{-\parsep}}
\newcommand{\ignore}[1]{}
\title{On linear programming bounds for spherical
codes and designs}
\author{Alex Samorodnitsky}

\begin{document}
\date{}
\maketitle
 
 
\newtheorem{THEOREM}{Theorem}[section]
\newenvironment{theorem}{\begin{THEOREM} \hspace{-.85em} {\bf :} 
}%
                        {\end{THEOREM}}
\newtheorem{LEMMA}[THEOREM]{Lemma}
\newenvironment{lemma}{\begin{LEMMA} \hspace{-.85em} {\bf :} }%
                      {\end{LEMMA}}
\newtheorem{COROLLARY}[THEOREM]{Corollary}
\newenvironment{corollary}{\begin{COROLLARY} \hspace{-.85em} {\bf 
:} }%
                          {\end{COROLLARY}}
\newtheorem{PROPOSITION}[THEOREM]{Proposition}
\newenvironment{proposition}{\begin{PROPOSITION} \hspace{-.85em} 
{\bf :} }%
                            {\end{PROPOSITION}}
\newtheorem{DEFINITION}[THEOREM]{Definition}
\newenvironment{definition}{\begin{DEFINITION} \hspace{-.85em} {\bf 
:} \rm}%
                            {\end{DEFINITION}}
\newtheorem{EXAMPLE}[THEOREM]{Example}
\newenvironment{example}{\begin{EXAMPLE} \hspace{-.85em} {\bf :} 
\rm}%
                            {\end{EXAMPLE}}
\newtheorem{CONJECTURE}[THEOREM]{Conjecture}
\newenvironment{conjecture}{\begin{CONJECTURE} \hspace{-.85em} 
{\bf :} \rm}%
                            {\end{CONJECTURE}}
\newtheorem{MAINCONJECTURE}[THEOREM]{Main Conjecture}
\newenvironment{mainconjecture}{\begin{MAINCONJECTURE} \hspace{-.85em} 
{\bf :} \rm}%
                            {\end{MAINCONJECTURE}}
\newtheorem{PROBLEM}[THEOREM]{Problem}
\newenvironment{problem}{\begin{PROBLEM} \hspace{-.85em} {\bf :} 
\rm}%
                            {\end{PROBLEM}}
\newtheorem{QUESTION}[THEOREM]{Question}
\newenvironment{question}{\begin{QUESTION} \hspace{-.85em} {\bf :} 
\rm}%
                            {\end{QUESTION}}
\newtheorem{REMARK}[THEOREM]{Remark}
\newenvironment{remark}{\begin{REMARK} \hspace{-.85em} {\bf :} 
\rm}%
                            {\end{REMARK}}
 
\newcommand{\thm}{\begin{theorem}}
\newcommand{\lem}{\begin{lemma}}
\newcommand{\pro}{\begin{proposition}}
\newcommand{\dfn}{\begin{definition}}
\newcommand{\rem}{\begin{remark}}
\newcommand{\xam}{\begin{example}}
\newcommand{\cnj}{\begin{conjecture}}
\newcommand{\mcnj}{\begin{mainconjecture}}
\newcommand{\prb}{\begin{problem}}
\newcommand{\que}{\begin{question}}
\newcommand{\cor}{\begin{corollary}}
\newcommand{\prf}{\noindent{\bf Proof:} }
\newcommand{\ethm}{\end{theorem}}
\newcommand{\elem}{\end{lemma}}
\newcommand{\epro}{\end{proposition}}
\newcommand{\edfn}{\bbox\end{definition}}
\newcommand{\erem}{\bbox\end{remark}}
\newcommand{\exam}{\bbox\end{example}}
\newcommand{\ecnj}{\bbox\end{conjecture}}
\newcommand{\emcnj}{\bbox\end{mainconjecture}}
\newcommand{\eprb}{\bbox\end{problem}}
\newcommand{\eque}{\bbox\end{question}}
\newcommand{\ecor}{\end{corollary}}
\newcommand{\eprf}{\bbox}
\newcommand{\beqn}{\begin{equation}}
\newcommand{\eeqn}{\end{equation}}
\newcommand{\wbox}{\mbox{$\sqcap$\llap{$\sqcup$}}}
\newcommand{\bbox}{\vrule height7pt width4pt depth1pt}
\newcommand{\qed}{\bbox}
\def\sup{^}

\def\H{\{0,1\}^n}

\def\S{S(n,w)}

\def\n{\lfloor \frac n2 \rfloor}

\def\Tp{Tchebyshef polynomial}
\def\Tps{TchebysDeto be the maximafine $A(n,d)$ l size of a code with distance $d$hef polynomials}
\newcommand{\rarrow}{\rightarrow}

\newcommand{\larrow}{\leftarrow}

\overfullrule=0pt
\def\setof#1{\lbrace #1 \rbrace}
 
\begin{abstract}
We investigate universal bounds on spherical codes and spherical
designs that could be obtained using Delsarte's 
linear programming methods. We
give a {\it lower} estimate for the LP upper bound on codes, and an  
{\it upper} estimate for the LP lower bound on designs.

Specifically, when the distance of the code is fixed and the
dimension goes to infinity, the LP upper bound on codes is at least as
large as the average of the best known upper and lower bounds.

When the dimension $n$ of the design is fixed, and the strength $k$
goes to infinity, the LP bound on designs turns out, 
in conjunction with known lower bounds, to be proportional to $k^{n-1}$.
%
\end{abstract}

\section{Introduction}
An $n$-dimensional spherical code of (angular) distance $\theta$ is a
subset of the $(n-1)$-dimensional unit sphere, such that the
angle between any two distinct points is at least
$\theta$. Equivalently, the Euclidean distance between any two
distinct points is at least  $2 \sin (\theta/2)$. 

An $n$-dimensional spherical design of strength $k$ is a
finite subset $W$ of the $(n-1)$-dimensional unit sphere, such that
for any algebraic polynomial $f$ of $n$ variables and degree $k$ holds
$$
\int_{S^{n-1}} f(x) dx = \frac{1}{|W|} \sum_{u\in W} f(u).
$$
We are interested
in the maximal cardinality $M(n,\theta)$ of a spherical code of
distance $\theta$, and in the minimal cardinality $N(n,k)$ of 
a design of strength $k$.

For a fixed $\theta$ and $n\rarrow \infty$, $M(n,\theta)$ increases
exponentially in $n$. 
The best known existential (lower) bound on the exponent
$\frac{1}{n} \log M(n,\theta)$ is
obtained by a volume 
argument \cite{CS}:
$$
\frac{1}{n} \log M(n,\theta) \ge \log \frac{1}{\sin \theta} - o(1),
$$
as $n$ goes to infinity. 

For a fixed $n$ and $k\rarrow \infty$, $M(n,k)$ increases
polynomially in $k$.
The best known existential (upper) bound on $N(n,k)$
is \cite{DT}
$$
N(n,k) \le O\left(k^{\frac{n(n-1)}{2}}\right), 
$$
as $k$ goes to infinity. 

The best universal bounds on codes and designs (upper for codes and
lower for designs) are obtained using linear programming methods,
initiated by Delsarte \cite{Dels}.

Let $\{P^{\alpha,\beta}_s\}$ be the Jacobi polynomials, orthogonal with
respect to a weight function $w^{\alpha,\beta}(t) =
(1-t)^{\alpha}(1+t)^{\beta}$ on 
$(-1,1)$. For $\alpha = \beta = \frac{n-3}{2}$, we will simply write
$\{P_s\}, ~w(t)$.	 
We will assume the standard normalization \cite{szego}, in
particular $P_0 \equiv 1$.  
Then \cite{KL}, \cite{DGS}
\beqn
\label{code_LP}
M(n,\theta) \le \min\left\{F(1):~F = \sum_{s=1}^{m}  a_s P_s;~ a_s \ge
0, ~a_0 = 1;~F\le 0 \mbox{ on } [-1, \cos \theta]\right\}
\eeqn
And \cite{DGS}
\beqn
\label{design_LP}
N(n,k) \ge \max\left\{\frac{1}{a_0}:~F = \sum_{s=1}^{m}  a_s P_s;~
F\ge 0 \mbox{ on } [-1, 1], ~F(1) = 1;~a_s \le 0 \mbox{ for } s\ge
k\right\}.
\eeqn
In (\ref{code_LP}) and (\ref{design_LP}) the degree $m$ of the
polynomial $F$ may be arbitrarily large.

We will denote the RHS of (\ref{code_LP}) by $M_{LP}(n,\theta)$ and 
the RHS of (\ref{design_LP}) by $N_{LP}(n,t)$. 

Kabatyansky and Levenshtein \cite{KL} obtain the best known upper
bound on $M(n,\theta)$ 
$$
\frac{1}{n} \log M_{LP}(n,\theta) \le \frac{1 + \sin{\theta}}{2 \sin
\theta} \log 
\frac{1 + \sin \theta}{2 
\sin \theta} - \frac{1 - \sin{\theta}}{2 \sin \theta} \log \frac{1 -
\sin{\theta}}{2 \sin \theta} + o(1).
$$
Yudin \cite{Y} gives the best known lower bound on $N(n,k)$, for $n$
fixed and $k\rarrow \infty$ 
$$
N_{LP}(n,k) \ge \frac{\int_{-1}^1 w(t)dt}{\int_{\gamma}^1 w(t)dt},
$$
where $\gamma$ is the maximal root of $P^{\frac{n-1}{2},\frac{n-1}{2}}_{k-1}$.
For $n$ fixed and $k\rarrow \infty$, this is at least (\cite{lev},
pp. 117-120) 
$$
\Omega\left(2^{-c (4n)^{\frac13}} \cdot \left(\frac{2}{n}\right)^n \cdot
k^{n-1}\right), 
$$
where $c \approx 1.86$, improving the lower
bound of Delsarte, Goethals, and Seidel \cite{DGS} 
by a factor of $\left(\frac{4}{e}\right)^n \cdot 
2^{-O\left(n^{\frac13}\right)}$.  

The exact values of $M_{LP}(n,\theta)$ and $N_{LP}(n,k)$ are not
known, and the relation of these derived quantities to $M(n,\theta)$
and $N(n,k)$ makes them legitimate subjects of research. In this
paper we obtain a lower bound on $M_{LP}(n,\theta)$ and an   
upper bound on $N_{LP}(n,k)$. This sets limits on how good the
bounds on codes and designs obtained through linear programming
methods could be. We follow the approach in 
\cite{sam}.

We prove:
\pro
\label{code}
For $n\ge 7$ holds 
\footnote{No significant attempt has been made, either here or in
proposition~\ref{design}, to extend the claim to the cases $n = 3,4,5,6$.}
\beqn
\label{LP_code_ineq}
M_{LP}(n,\theta) \ge 
\Omega\left(\frac{1}{rn^{1/2}}\right)
\left(\frac{1}{1-\delta^2}\right)^{\frac{n-4}{4}} \frac{P_r(1)}{\|P_r\|_2},  
\eeqn
where
$\delta = \cos \theta$, 
$r:=\max\left\{s~:~x_s \le \delta - 2n^{-\frac12}\right\}$, and 
$x_s$ denotes the maximal root of $P_s$.
\epro
\pro
\label{design}
Let $\ell = k$ if $k$ is even, and $\ell = k+1$ if $k$ is odd. Then
for $n \ge 6$ holds
\beqn
\label{LP_design_ineq_root}
N_{LP}(n,k) \le O\left(k\right) 
\left(\frac{1}{1-\rho^2}\right)^{\frac{n-2}{4}} \cdot
\frac{P_{\ell}(1)}{\|P_{\ell}\|_2},   
\eeqn
where $\rho$ is the maximal root of $P^{\frac{n-5}{2},\frac{n-5}{2}}_{\ell}$.
\footnote{Observe that bounds (\ref{LP_code_ineq}) and
(\ref{LP_design_ineq_root}) are, in a certain sense, up
to polynomial factors, 'dual' to
each other, as are linear programs (\ref{code_LP}) and
(\ref{design_LP}).} 
\epro

Analyzing the asymptotic behaviour of the bounds leads to
following corollaries.
\cor
\label{code_asymptotic}
$$
\frac{1}{n} \log M_{LP}(n,\theta) \ge \frac{\log \frac{1}{\sin \theta}}{2} + 
\frac{\frac{1 + \sin{\theta}}{2 \sin \theta} \log \frac{1 + \sin \theta}{2
\sin \theta} - \frac{1 - \sin{\theta}}{2 \sin \theta} \log \frac{1 -
\sin{\theta}}{2 \sin \theta} }{2} + o(1),
$$
as $n$ goes to infinity.
\ecor
\cor
\label{design_asymptotic}
$$
N_{LP}(n,k) \le O\left(n^{-\frac12} \cdot
\left(\frac{\sqrt{2e}}{n}\right)^{n-3} \cdot k^{n-1}\right),
$$
for $n$ fixed and $k$ going to infinity.
\ecor
Combining this upper bound with Yudin's lower bound on
$N_{LP}(n,k)$ (or with the lower bound of Delsarte,
Goethals, and Seidel) we obtain
$$
N_{LP}(n,k) = \Theta\left(k^{n-1}\right),
$$
for $n$ fixed and $k$ going to infinity.

Note that the bound in (\ref{code_asymptotic}) is,
asymptotically, a geometric mean of the existential lower
bound and the  Kabatyansky-Levenshtein upper bound on $M(n,\theta)$.
Similar situation \cite{sam}
holds in the context of the LP-bounds for binary and constant weight
binary codes. 

The paper is organized as follows: in the next section we provide
relevant information about Jacobi polynomials. Propositions \ref{code}
and \ref{design} are proved in sections \ref{sec_code} and \ref{sec_design}.

\section{Preliminaries}
We will require some facts about Jacobi polynomials
$P^{\alpha, \alpha}_s$. These facts are presented in this section.

\noindent {\bf Normalization} \cite{szego}\\
\beqn
\label{value_at_1}
P_s(1) = \left(\frac{n-1}{2}\right)_s,
\eeqn
where $(x)_s := x(x+1)...(x+s-1)$.
\beqn
\label{norm}
\|P_s\|^2_2 = \frac{2^{n-2}}{2s + n - 2} \frac{\Gamma^2\left(r +
\frac{n-1}{2}\right)}{r! (r+ n-3)!}.
\eeqn
\noindent {\bf Asymptotics of the maximal root} \cite{lev}\\
Let $x_s$ be the maximal root of $P_s$. Then estimates in \cite{lev}
(cor. 5.17, identity 5.35) give:
$$
\left(2\sqrt{s} - 1\right) \sqrt{\frac{s+n-4}{(2s+n-6)(2s+n-4)}}
\le x_s \le 2\sqrt{s-1} \sqrt{\frac{s+n-4}{(2s+n-6)(2s+n-4)}}.
$$
It follows that for any $s>0$, and $n\ge 6$
\beqn
\label{root_estimate}
\bigg | x_s - \frac{\sqrt{4s(s+n)}}{2s+n} \bigg | \le
\sqrt{\frac 2n}.
\eeqn
It also follows that for any $s > 1$ and $n \ge 4$, 
\beqn
\label{root_upperbound}
1 - x^2_s \ge \frac{(n-4)^2}{(2s + n - 4)^2}.
\eeqn

\lem
Assuming $n\ge 6$, for any $s \ge 0$ holds: 
$w(t) P^2_s(t)$ is a decreasing function of $t$ in the interval 
$\left[x_s+ 2n^{-\frac12}, 1\right]$.  
\elem
\prf
Let $t \in \left[x_s+ 2n^{-\frac12}, 1\right]$. By
(\ref{root_estimate}), $t \ge 
\frac{\sqrt{4s(s+n)}}{2s+n} + \frac12 n^{-\frac12}$. It is not hard to
check that this implies $(n-1)^2 t^2 - 4(1-t^2)s(s+n-2) \ge 4t^2$. 

Now we can follow the analysis of \cite{ABL} for
$\frac{P'_s(t)}{P_s(t)}$, obtaining:
$$
\frac{P'_s(t)}{P_s(t)} <  \frac{(n-1)t - \sqrt{(n-1)^2 t^2 -
4(1-t^2)s(s+n-2)}}{2(1-t^2)} < \frac{(n-3)t}{2(1-t^2)}.
$$
We conclude the proof of the lemma by computing $\frac{d}{dt} \ln(w(t)
P^2_s(t)) = 2\frac{P'_s(t)}{P_s(t)} - \frac{(n-3)t}{1-t^2} <
0$, for $t \in \left[x_s + 2n^{-\frac12}, 1\right]$.
\eprf

\cor
\label{monotone}
Assuming $n\ge 6$, for any $s \ge 0$ holds: 
$w(t) P_s(t)$ is a decreasing function of $t$ in the interval 
$\left[x_s + 2n^{-\frac12}, 1\right]$. 
\ecor 

\lem
\label{first}
Assuming $n \ge 7$, for any $r > 0$ holds:
$$
\sum_{s = r}^{\infty} \frac{\|P_s\|_2}{P_s(1)} \le O(r) \cdot 
\frac{\|P_r\|_2}{P_r(1)}.
$$
\elem
\prf
We will assume $n$ is odd, the proof for even $n$ is similar.
Set $a_s = \frac{\|P_s\|}{P_s(1)}$. 
By \cite{szego} 
$$
a^2_s = \frac{\frac{2^{n-2}}{2s+n-2} \cdot 
\frac{(s + (n-3)/2)!(s + (n-3)/2)!}{s! (s+n-3)!}}{{{s +
(n-3)/2}\choose s}^2}  = 
\frac{2^{n-2}}{{{n-3}\choose (n-3)/2}} \cdot \frac{1}{2s+n-2} \cdot 
\frac{1}{{{s+n-3}\choose s}}.
$$
Therefore $\frac{a^2_{s+1}}{a^2_s} = \frac{2s+n-2}{2s+n}\cdot
\frac{s+1}{s+n-2} \le \frac{s+1}{s+n-2}$, and for any $t\ge 0$, 
$$
\frac{a^2_{s+t}}{a^2_s} = \prod_{i=0}^{t-1}
\frac{a^2_{s+i+1}}{a^2_{s+i}} \le
\frac{(s+1)_t}{(s+n-2)_t} = 
\frac{(s+1)_{n-3}}{(s+t+1)_{n-3}}.
$$
It follows that
$$
\sum_{s = r}^{\infty} a_s \le a_r \cdot \sum_{t=0}^{\infty}
\sqrt{\frac{(r+1)_{n-3}}{(r+t+1)_{n-3}}}.
$$
The last sum, assuming $n \ge 7$, is at most $a_r \cdot \left[ O(r) + O(r^2)
\sum_{k=r}^{\infty} \frac{1}{k^2}\right] = O(r) a_r$. 
\eprf 

\rem
\label{ratio_monoton}
Observe that the ratio $\frac{\|P_s\|_2}{P_s(1)}$
decreases with $s$.
\erem

\thm \cite{NEM}
For all $x\in [-1,1]$
\beqn
\label{maximum}
|P_s(x)| \le O(\sqrt{n})
\frac{\|P_s\|_2}{\left(1-x^2\right)^{\frac{n-2}{4}}}.
\eeqn
\ethm

\lem
\label{design_monotone}
Assuming $n\ge 3$, for any $t > 0$ holds: 
$\frac{P_s(1) P_s(t)}{\|P_s\|^2_2}$ is increasing in $s$ for even $s$ 
such that the
maximal root of the Jacobi polynomial
$P^{\frac{n-5}{2},\frac{n-5}{2}}_{s+2}$ does not exceed $t$. 
\elem
\prf
Let $b_s =  \frac{P_s(1) P_s(t)}{\|P_s\|^2_2}$. We have to prove that
$b_{s+2} \ge b_s$, which is equivalent to 
\beqn
\frac{P_{s+2}(1)}{P_s(1)}
\cdot \frac{\|P_s\|^2_2}{\|P_{s+2}\|^2_2} \cdot P_{s+2}(t) \ge
P_s(t).
\label{jacobi_ineq}
\eeqn
It will be useful to renormalize and work with the ultraspherical polynomials
$C_s = C^{(\frac{n-2}{2})}_s$, which are proportional to Jacobi
polynomials $P^{\frac{n-3}{2},\frac{n-3}{2}}_s$: 
$$
C_s = \frac{\Gamma\left(\frac{n-1}{2}\right) \Gamma(s+ n -2)}
{\Gamma(n-2) \Gamma\left(s + \frac{n-1}{2}\right)} \cdot P_s
$$
Rewriting (\ref{jacobi_ineq}) for ultraspherical polynomials, and
substituting the values of $P_i(1)$ and $\|P_i\|^2_2$, for $i = s,
s+2$, we get the
following inequality to prove:
$$
(2s + n +2) C_{s+2}(t) \ge (2s + n -2) C_s(t).
$$
Consider the following identity (\cite{bateman},
p. 178, (36)):
$$
\frac{n-4}{2} \cdot (C_{s+2}(t) - C_s(t)) = \left(s +
\frac{n-2}{2}\right) C^{(\frac{n-4}{2})}_{s+2}(t).
$$
In the assumed range for $s$, $C^{(\frac{n-4}{2})}_{s+2}(t) \ge 0$.
Therefore $C_{s+2}(t) \ge C_s(t)$. In order to complete the proof it
is sufficient to show that $C_{s+2}(t) = C^{(\frac{n-2}{2})}_s(t) \ge
0$. This is indeed true, because by a theorem of Markov (\cite{szego},
(6.21.3)), if $\lambda > \beta$, then the maximal root of
$C^{(\lambda)}_s$ is smaller than 
that of $C^{(\beta)}_s$. 
\eprf

\section{A lower bound on $M_{LP}(n,\theta)$}
\label{sec_code}
\prf {\bf of proposition~\ref{code}}\\
First, 
$$
\int_{\delta}^1 F(t) w(t) dt \ge  
\int_{-1}^1  F(t) w(t) dt = a_0 \int_{-1}^1  P^2_0(t) w(t) dt = 
\int_{-1}^1 w(t) dt.
$$
Therefore, for some $t_0 \in [\delta, 1]$ holds
$F(t_0) w(t_0) \ge \frac{1}{1-\delta} \int_{-1}^1 w(t)
dt$. 
Let $F = F_1 + F_2$, where $F_1 := \sum_{s=0}^r a_s P_s$. 

We would like to show that either 
$$
F_2(t_0) w(t_0) \ge \frac{1}{2(1-\delta)} \int_{-1}^1 w(t) dt
$$
or 
$$
\Big | F_2(\delta) w(\delta) \Big | \ge \frac{1}{2(1-\delta)}
\int_{-1}^1 w(t) dt. 
$$ 
If 
$F_1(t_0) w(t_0) \le \frac{1}{2(1-\delta)} \int_{-1}^1 w(t)
dt$, then the first inequality holds. Otherwise, by 
corollary~\ref{monotone},
$$
F_1(\delta) w(\delta) \ge F_1(t_0) w(t_0) > \frac{1}{2(1-\delta)}
\int_{-1}^1 w(t) dt. 
$$ 
Since $F(\delta) \le 0$, it must be that 
$F_2(\delta) w(\delta)  < - \frac{1}{2(1-\delta)} \int_{-1}^1 w(t) dt$,
and the second inequality holds.

Let $t_m$ be one of the two points $t_0$, $\delta$, so that $\big |
F_2(t_m) w(t_m) \big | \ge  
\frac{1}{2(1-\delta)} \int_{-1}^1 w(t) dt$.
Then, using (\ref{maximum})
$$
\frac{1}{2(1-\delta)} \int_{-1}^1 w(t) dt \le \Big | F_2(t_m) w(t_m)
\Big | \le w(t_m) \cdot \sum_{s=r+1}^{m} a_s  \Big |P_s(t_m)\Big |\le 
$$
$$
(1-t^2_m)^{\frac{n-4}{4}} \cdot \sum_{s=r+1}^{m} a_s \|P_s\|_2 \le
(1-\delta^2)^{\frac{n-4}{4}} \cdot 
\sum_{s=r+1}^{m} a_s \|P_s\|_2.
$$
Since all the coefficients $a_s$ are nonnegative, they are bounded
from above: $a_s \le \frac{F(1)}{P_s(1)}$. Therefore
$$
F(1) \ge \frac{1}{2(1-\delta)} \int_{-1}^1 w(t) dt \cdot
\left(\frac{1}{1-\delta^2}\right)^{\frac{n-4}{4}}  \cdot
\frac{1}{\sum_{s=r+1}^{m} \frac{\|P_s\|_2}{P_s(1)}} \ge 
\Omega\left(\frac{1}{rn^{1/2}}\right)
\left(\frac{1}{1-\delta^2}\right)^{\frac{n-2}{4}} \frac{P_r(1)}{\|P_r\|_2}.  
$$
The last inequality uses lemma~\ref{first} and a simple fact: 
$\int_{-1}^1 w(t) dt = \Theta\left(n^{-\frac12}\right)$.
\eprf

\prf {\bf of corollary~\ref{code_asymptotic}}\\
The main step is to estimate $r$. From (\ref{root_estimate}), 
$\Big | r - \frac{\frac{1}{1-\delta^2} - 1}{2} \cdot n \Big
| \le O\left(n^{\frac12}\right)$. 
Now, the claim of the corollary is obtained using (\ref{value_at_1})
and (\ref{norm}), and simplifying.
\eprf 

\section{An upper bound on $N_{LP}(n,k)$}
\label{sec_design}
\prf {\bf of proposition~\ref{design}}\\
First, we may, without loss of generality, 
assume that $F$ is symmetric around
zero. Indeed, if $F$ is not symmetric, consider the symmetric function $G =
\frac{F(t) + F(-t)}{F(1) + F(-1)}$. Clearly
$G\ge 0$ on $[-1, 1]$, $G(1) = 1$, and in the expansion 
$G = \sum_{s=0}^{\infty}  b_s P_s$, the coefficients $b_s$ are
nonpositive for $s \ge k$. Also $b_0 = \frac{2 a_0}{F(1) + F(-1)} \le
2 a_0$, so it is sufficient to provide a lower bound for $b_0$.

This said, we assume that initially $F$ is symmetric. This, in particular,
implies that $a_s = 0$ for all odd $s$.

To make this proof as similar as possible, up to a
'duality', to the proof of proposition~\ref{code}, we introduce
two definitions: Let $A : {\bf N} \rarrow {\bf R}$ be defined by $A(s)
= a_s$. Then 
$$
A(s) = \frac{\int_{-1}^1 F(t) P_s(t) w(t) dt}{\|P_s\|^2_2} =
\int_{-1}^1 F(t) \alpha_t(s) dt, 
$$
where 
$$
\alpha_t(s) = \frac{P_s(t) w(t)}{\|P_s\|^2_2}
$$
is 'dual' to $P_s(t)$. 

Now, $\sum_{s=0}^{k-1} P_s(1)
A(s) \ge \sum_{s=0}^{m} P_s(1) A(s) = F(1) = 1$. Therefore, 
there exists an index $s_0 \in [0,k-1]$ such that
$P_{s_0}(1) A(s_0) \ge \frac{1}{k}$.

Write $A = A_1 + A_2$, where $A_2(s) := \int_{-\rho}^{\rho} F(t)
\alpha_t(s) dt$. Let $\ell = k$ if $k$ is even, and $\ell = k+1$ if
$k$ is odd. We would like to show that either
$$
P_{s_0}(1) A_2(s_0) \ge \frac{1}{2k},
$$
or 
$$
P_{\ell}(1) \Big | A_2(\ell) \Big | \ge \frac{1}{2k}.
$$
If 
$P_{s_0}(1) A_1(s_0) \le \frac{1}{2k}$, then the first inequality
holds. Otherwise, observe that, by lemma~\ref{design_monotone}, for
every $t \in 
[\rho,1]$ holds $P_{\ell}(1) \alpha_t(\ell) \ge P_{s_0}(1)
\alpha_t(s_0)$. 
This is also true for all $t \in [-1,-\rho]$, since $s_0$, $\ell$ are even
and consequently $P_{s_0}, P_{\ell}$ are symmetric around $0$. 
Therefore
$$
P_{\ell}(1) A_1(\ell) = P_{\ell}(1) \int_{[-1,\rho] \cup [\rho,1]} F(t)
\alpha_t(\ell) dt \ge 
P_{s_0}(1) \int_{[-1,\rho] \cup [\rho,1]} F(t) \alpha_t(s_0) dt >
\frac{1}{2k}.
$$
Since $A(\ell) \le 0$, it must be that $P_{\ell}(1) \big | A_2(\ell) \big
| \ge P_{\ell}(1) 
A_1(\ell) > \frac{1}{2k}$, and the second inequality holds. 

Let $s$ be one of the two indices $s_0, \ell$, so that $P_s(1) \big
| A_2(s) \big | \ge \frac{1}{2k}$. Then
$$
\frac{1}{2k} \le P_s(1) \Big | A_2(s) \Big | = P_s(1) \cdot \bigg |
\int_{-\rho}^{\rho} F(t) \alpha_t(s) dt \bigg | = \frac{P_s(1)}{\|P_s\|^2_2} 
\cdot \bigg | \int_{-\rho}^{\rho} F(t) P_s(t) w(t) dt \bigg | 
$$
$$
\le 
\frac{P_s(1)}{\|P_s\|^2_2} \cdot \max_{t \in [-\rho,\rho]} \Big | P_s(t)
\Big |  \int_{-1}^1 F(t) w(t) dt = A(0) \cdot \int_{-1}^1 w(t) dt  
\cdot \frac{P_s(1)}{\|P_s\|^2_2} \cdot \max_{t \in [-\rho,\rho]} \Big | P_s(t)
\Big | \le 
$$
$$
A(0) \cdot \Omega(1) \cdot
\left(\frac{1}{1-\rho^2}\right)^{\frac{n-2}{4}} \cdot
\frac{P_s(1)}{\|P_s\|_2}.  
$$
The last inequality uses (\ref{maximum}).

Therefore,
$$
\frac{1}{a_0} = \frac{1}{A(0)} \le
O\left(k\right) \cdot \left(\frac{1}{1-\rho^2}\right)^{\frac{n-2}{4}}
\cdot \frac{P_s(1)}{\|P_s\|_2} \le
O\left(k\right) \cdot \left(\frac{1}{1-\rho^2}\right)^{\frac{n-2}{4}}
\cdot \frac{P_{\ell}(1)}{\|P_{\ell}\|_2}, 
$$
using the fact that the fraction
$\frac{P_s(1)}{\|P_s\|_2}$ is increasing in $s$ (see
remark~\ref{ratio_monoton}).  
\eprf

\prf {\bf of corollary~\ref{design_asymptotic}}\\
The claim of the corollary is obtained using
(\ref{value_at_1}), (\ref{norm}), and (\ref{root_upperbound}), and
simplifying.    
\eprf


\begin{thebibliography}{99}

\bibitem{ABL}
A. Ashikhmin, A. Barg, and S. N. Litsyn, {\sl A new upper bound on the
reliability function of the Gaussian channel}, IEEE
Trans. Inform. Theory, vol. IT-46, 2000, 1945-1961.



\bibitem{bateman}
H. Bateman, A. Erd\'elyi, {\bf Higher transcedental functions}, Vol. 2,
McGraw-Hill, New York, 1953. 

\bibitem{BE}
P. Borwein and T. Erd\'elyi, {\bf Polynomials and Polynomial Inequalities},
Springer-Verlag, New York, 1995.

\bibitem{CS}
J. H. Conway and N. J. A. Sloane, {\bf Sphere Packings, Lattices and
Groups}, Springer-Verlag, New York, 1999.

\bibitem{Dels}
P. Delsarte, {\sl An algebraic approach to the association schemes of
coding theory}, Philips Res. Rep. Suppl. No. 10 (1973)

\bibitem{DGS}
P. Delsarte, J. M. Goethals and J. J. Seidel, {\sl Spherical codes and
designs}, Geom. Dedic., 6, 1979, 363-388.

\bibitem{DT}
M. Drmota and R. F. Tichy, {\bf Sequences, Discrepancies and
Applications}, Lecture Notes
in Mathematics, 1651. Springer-Verlag, Berlin, 1997.

\bibitem{KL}
G. Kabatyansky and V. I. Levenshtein, {\sl Bounds for packings on the
sphere and in the space}, Probl. Pered. Inform., vol. 14, no. 1, 1978,
3-25.

\bibitem{lev}
V. I. Levenshtein, {\sl Universal bounds for codes and designs}, in
Handbook of Coding Theory, V.S. Pless 
and W.C. Huffman Eds., Amsterdam, Elsevier, 1998. 

\bibitem{NEM}
P. Nevai, T. Erd\'elyi, and A. P. Magnus, {\sl Generalized Jacobi
weights, Christoffel functions, and Jacobi polynomials}, SIAM
J. Math. Anal., vol. 25, no. 2, 1994, 602-614.

\bibitem{sam} 
A. Samorodnitsky, {\sl On the optimum of Delsarte's liner program},
J. Comb. Th., Ser. A, 96, no. 2, 2001, 261-287. 

\bibitem{szego}
G. Szeg\"o, {\bf Orthogonal Polynomials}, American Mathematical Society, 1939.

\bibitem{Y}
V. A. Yudin, {\sl Lower bounds for spherical designs},
Izvestiya. Mathematics 61:3, 1997, 673-683.

\end{thebibliography}
\end{document}